\documentclass{cmslatex}
\usepackage[utf8]{inputenc}
\usepackage[paperwidth=7in, paperheight=10in, margin=.875in]{geometry}
\usepackage[backref,colorlinks,linkcolor=red,anchorcolor=green,citecolor=blue]{hyperref}
\usepackage{booktabs,multirow}
\usepackage[noend]{algpseudocode}
\usepackage[pagewise]{lineno}
\usepackage{algorithmicx,algorithm,setspace,subfigure}
\usepackage{amsfonts,amssymb}
\usepackage{amsmath}
\usepackage{graphicx}
\usepackage{amsmath}
\usepackage{cite}
\usepackage{subfigure}
\usepackage{todonotes}
\usepackage{enumerate}
\usepackage{ulem,marvosym}
\renewcommand{\theequation}{\arabic{section}.\arabic{equation}}

\newcommand{\abs}[1]{\left|#1\right|}
\newcommand{\brac}[1]{\left(#1\right)}
\newcommand{\norm}[1]{\left\Vert#1\right\Vert}

\newcommand{\bp}{\boldsymbol{p}}
\newcommand{\bq}{\boldsymbol{q}}
\newcommand{\br}{\boldsymbol{r}}
\newcommand{\bs}{\boldsymbol{s}}
\newcommand{\mbr}{\mathbb{R}}
\newcommand{\bu}{\boldsymbol{u}}
\newcommand{\bv}{\boldsymbol{v}}
\newcommand{\bgamma}{\boldsymbol{\gamma}}
\newcommand{\balpha}{\boldsymbol{\alpha}}
\newcommand{\bbeta}{\boldsymbol{\beta}}
\newcommand{\bphi}{\boldsymbol{\phi}}
\newcommand{\bpsi}{\boldsymbol{\psi}}
\newcommand{\veps}{\varepsilon}
\allowdisplaybreaks
\begin{document}
\begin{sloppypar}

 \title{Fast Sinkhorn I: an $O(N)$ algorithm for the Wasserstein-1 metric \thanks{Received date, and accepted date (The correct dates will be entered by the editor).}}


          \author{Qichen Liao\thanks{Department of Mathematical Sciences, Tsinghua University, Beijing 100084, China; Theory Lab, Central Research Institute, 2012 Labs, Huawei Technologies Co. Ltd., Hong Kong SAR, China (\href{mailto:lqc20@mails.tsinghua.edu.cn}{lqc20@mails.tsinghua.edu.cn}).}
          \and Jing Chen\thanks{School of Physical and Mathematical Sciences, Nanyang Technological University, Singapore 639798 (\href{mailto:jing.chen@ntu.edu.sg}{jing.chen@ntu.edu.sg}).}
          \and Zihao Wang\thanks{Department of Computer Science and Engineering, Hong Kong University of Science and Technology, Clear Water Bay, Hong Kong SAR, China (\href{mailto:zwanggc@cse.ust.hk}{zwanggc@cse.ust.hk}).}
          \and Bo Bai\thanks{Theory Lab, Central Research Institute, 2012 Labs, Huawei Technologies Co. Ltd., Hong Kong SAR, China (\href{mailto:baibo8@huawei.com}{baibo8@huawei.com}).}
          \and Shi Jin\thanks{Corresponding author. School of Mathematical Sciences, Institute of Natural Sciences, and MOE-LSC Shanghai Jiao Tong University, Shanghai 200240, China (\href{mailto:shijin-m@sjtu.edu.cn}{shijin-m@sjtu.edu.cn}).}
          \and Hao Wu\thanks{Corresponding author. Department of Mathematical Sciences, Tsinghua University, Beijing 100084, China (\href{mailto:hwu@tsinghua.edu.cn}{hwu@tsinghua.edu.cn})}.}

         \pagestyle{myheadings} \markboth{THE FS-1 ALGORITHM FOR THE WASSERSTEIN-1 METRIC}{Q. LIAO, J. CHEN, Z. WANG, B. BAI, S. JIN AND H. WU} \maketitle
\begin{abstract}
    The Wasserstein metric is broadly used in  optimal transport   for comparing two probabilistic distributions, with successful applications in various fields such as machine learning, signal processing, seismic inversion, etc. Nevertheless, the high computational complexity is an obstacle for  its practical applications. The Sinkhorn algorithm, one of the main methods in computing the Wasserstein metric,  solves an entropy regularized minimizing problem, which allows arbitrary approximations to the Wasserstein metric with $O(N^2)$ computational cost. However, higher accuracy of its numerical approximation requires more Sinkhorn iterations with repeated matrix-vector multiplications, which is still unaffordable. In this work, we propose an efficient implementation of the  Sinkhorn algorithm to calculate the Wasserstein-1 metric with $O(N)$ computational cost, which achieves the optimal theoretical complexity. By utilizing the special structure of Sinkhorn's kernel, the repeated matrix-vector multiplications can be implemented with $O(N)$ times multiplications and additions, using the Qin Jiushao or Horner's method for efficient polynomial evaluation,  leading to an efficient algorithm without losing accuracy. In addition, the log-domain stabilization technique, used to stabilize the iterative procedure, can also be applied in this  algorithm. Our numerical experiments show that the newly developed  algorithm is one to three orders of magnitude faster than the original Sinkhorn algorithm.
\end{abstract}

\begin{keywords}  
	Optimal Transport, Wasserstein-1 metric, Sinkhorn algorithm, FS-1 algorithm, fast algorithm
\end{keywords}

\begin{AMS} 
 	49M25; 49M30; 65K10
\end{AMS}

\section{Introduction}

The Wasserstein metric has been widely used in optimal transport  for the global comparison between probabilistic distributions. It has  been successfully used  in various fields such as machine learning~\cite{goodfellow2014generative,lin2021wasserstein}, image processing~\cite{rubner2000earth}, inverse problems~\cite{chen2018quadratic,engquist2020quadratic,metivier2016measuring,yang2018application,heaton2020wasserstein}, and density function theory\cite{hu2021global,buttazzo2012optimal,cotar2013density}.
 Many numerical methods have been developed, including the linear programming methods~\cite{pele2009fast,li2020asymptotically,yang2021fast}, combinatorial methods~\cite{santambrogio2015optimal}, solving Monge-Amph\`ere equations~\cite{froese2011convergent,froese2012numerical,benamou2014numerical,benamou2000computational} and proximal splitting methods~\cite{combettes2011proximal,metivier2016optimal}. In recent years, several approximation techniques in optimal transport for high-dimensional distributions have also been proposed approximately~\cite{meng2019large,meng2020sufficient}.

One of the popular  numerical techniques to compute the Wasserstein metric, is the the Sinkhorn algorithm~\cite{cuturi2013sinkhorn,sinkhorn1967diagonal}, which minimizes the entropy regularized problem. It provides the solution roughly in $O(N^2)$ operations with guaranteed convergence~\cite{lin2019efficiency}.  With the help of GPU acceleration, the efficiency of using the Sinkhorn algorithm to solve the optimal transport problem can be more significantly improved~\cite{li2018parallel,ryu2018vector,ryu2018unbalanced}. For the Wasserstein-2 matric, the computation can be accelerated using the Gaussian convolution by approximating the geodesic distance-based kernel with the heat kernel~\cite{solomon2015convolutional}. Moreover, through statistical sampling, dimensional reduction, and other approximation methods, the complexity of the Sinkhorn algorithm can be reduced to $O(N\log N)$~\cite{altschuler2018massively,klicpera2021scalable} or even $O(N)$~\cite{scetbon2020linear}. An alternative way is to define the Wasserstein metric on the finite tree space, and the computational complexity could be $O(N\log N)$~\cite{DBLP:conf/nips/LeYFC19,DBLP:conf/icml/TakezawaSY21}. Because of these progresses, the Sinkhorn algorithm has been widely used in practical problems. However, those techniques provide just approximations of the transportation cost, rather than the precise computation of the original optimization problem. Therefore, it is still of great interest to develop fast and accurate Sinkhorn-type algorithms for solving large-scale optimal transport problems.

In recent years, some   fast algorithms for accurately solving the optimal transport problem have been developed recently. For example, through multi-level grids, the computation complexity can be reduced to $O(N^{1.5}\log N)$~\cite{liu2021multilevel}. Another well-known conclusion is that for the one-dimensional quadratic Wasserstein metric, the complexity of the sorting-based algorithm is only $O(N\log N)$~\cite{peyre2019computational}.

In this work, by observing the special structure of the kernel matrix in the Sinkhorn algorithm to solve the Wasserstein-1 metric, we propose a  novel matrix-vector multiplication based on dynamic programming~\cite{kleinberg2006algorithm}. During each iteration step, it involves only a forward and backward recursive sweeping process, as in Qin Jiushao's (or Horners', although it appeared several centuries later) method for polynomial evaluations \cite{needham1974, horner1819}, which  reduces the computational cost  from $O(N^2)$ to $O(N)$ in each step, thus  achieving the optimal theoretical complexity.  In addition, the log-stabilization~\cite{chizat2018scaling},  an important technique to improve the numerical stability of the Sinkhorn algorithm, can be implemented in this strategy with the same stability property. In this paper we abbreviate this method as FS-1. 

The rest of the paper is organized as follows. In Section \ref{sec:2}, the basics of the Wasserstein-1 metric and the Sinkhorn algorithm are briefly introduced. After showing the elaborate structure of the obtained kernel matrix, we use it to develop the FS-1 algorithm in Section \ref{sec:3}. We will also analyze the stability and integrate the log-stabilization technique into our FS-1 to improve the numerical stability. In Section \ref{sec:4}, the FS-1 is generalized to higher dimension. We will provide numerical experiments to illustrate the huge efficiency advantage of our FS-1 algorithm in Section \ref{sec:5}. Finally, we conclude the paper in Section \ref{sec:6}.

\section{The Wasserstein-1 metric and the Sinkhorn algorithm}
\label{sec:2}
    Consider two probabilistic density functions $u(x)$ and $v(x)$ on a doman $\Omega\subset \mbr$, the Kantorovich's formulation of Wasserstein-1 metric is defined as~\cite{villani2009optimal}:
	\begin{equation}
		\label{kantorovich}
		\begin{split}
		W(u,v) &= \inf_{\gamma(x,y) \in \Gamma} \int_{\Omega \times \Omega} d(x,y) \gamma(x,y) \mathrm{d}x \mathrm{d}y, \\
		\Gamma = \Big\{ \gamma(x,y) \Big| &\int_\Omega \gamma(x,y) \mathrm{d}y = u(x), \int_{\Omega} \gamma(x,y) \mathrm{d}x = v(y) \Big\},
		\end{split}
	\end{equation}
	where $d(x, y) = |x - y|$.
    The Sinkhorn algorithm proposed in~\cite{cuturi2013sinkhorn,sinkhorn1967diagonal} introduces an entropy term and solves the regularized problem:
	\begin{equation}\label{eqn:entropy problem}
	    	W_\veps(u,v) = \inf_{\gamma(x,y) \in \Gamma} \int_{\Omega \times \Omega} \abs{x-y} \gamma(x,y) + \veps \gamma(x,y) \ln\big( \gamma(x,y) \big) \mathrm{d}x \mathrm{d}y.
	\end{equation}
	As for numerical realization, we consider two discretized probabilistic distributions 
	\begin{equation*}
	\bu = (u_1,u_2,\cdots u_N), \quad  \bv = (v_1,v_2,\cdots, v_N),
	\end{equation*}
	on a uniform mesh grid with a grid spacing of $h$. Then the entropy regularized minimizing problem \eqref{eqn:entropy problem} can be discretized as 
	\begin{equation}\label{eqn:disc entropy problem}
		W_\veps(\bu,\bv) = \inf_{\gamma_{ij}} \sum_{i=1}^N \sum_{j=1}^N \brac{\gamma_{ij}\abs{i-j}h + \veps \gamma_{ij} \ln\big( \gamma_{ij} \big)},
 	\end{equation} 
	and $\gamma_{ij}$ satisfies	
	\begin{equation*}\label{gammaij}
		\sum_{j=1}^N \gamma_{ij} = u_i, \quad \sum_{i=1}^N \gamma_{ij} = v_j, \quad \gamma_{ij}\geq0.
	\end{equation*} 
    The Lagrangian of the above equations writes
    \begin{equation*}\label{Laplace}
        L(\bgamma,\balpha,\bbeta) = \sum_{i=1}^N \sum_{j=1}^N \brac{\gamma_{ij}\abs{i-j}h + \veps \gamma_{ij} \ln(\gamma_{ij})} + \sum_{i=1}^N \alpha_i (\sum_{j=1}^N \gamma_{ij}-u_i) + \sum_{j=1}^N \beta_j (\sum_{i=1}^N \gamma_{ij}-v_j).
    \end{equation*} 
    Taking the derivative of the Lagrangian with respect of $\gamma_{ij}$ directly leads to
    \begin{equation*}
        \gamma_{ij} = e^{-\frac{1}{2}-\frac{1}{\veps}\alpha_i} K_{ij} e^{-\frac{1}{2}-\frac{1}{\veps}\beta_j}, \;\;
        \textrm{where} \;\; K_{ij}=e^{-\abs{i-j}h/\veps}.
    \end{equation*}
    To avoid $\veps$ in the denominator, set $\phi_i=e^{-\frac{1}{2}-\frac{1}{\veps}\alpha_i}$ and $\psi_j=e^{-\frac{1}{2}-\frac{1}{\veps}\beta_j}$, we obtain
    \begin{equation*}\label{solution}
        \phi_i \sum_{j=1}^N K_{ij} \psi_j = u_i, \quad
        \psi_j \sum_{i=1}^N K_{ij} \phi_i = v_j, \quad \forall i,\;j=1,2,\cdots,N.
    \end{equation*}
    Since the entries in $K=(K_{ij})$ are strictly positive, Sinkhorn's iteration \cite{sinkhorn1967diagonal} can be applied to iteratively update vectors $\bphi=(\phi_i)$ and $\bpsi=(\psi_j)$ by pointwise computation:
    \begin{equation}\label{update_sinkhorn}
        \bpsi^{(\ell+1)} = \bv\oslash (K^T\bphi^{(\ell)}),\quad
        \bphi^{(\ell+1)} = \bu\oslash (K\bpsi^{(\ell+1)}),
    \end{equation}
    in which the notion $\oslash$ represents pointwise division and $\ell$ denotes the iterative steps.
    

\medskip

\begin{remark} \label{rmk_stb}
     In~\cite{chizat2018scaling}, a log-domain stabilization technique is proposed to reduce the numerical instability caused by the small parameter $\veps$. The idea is that when the infinite norms of $\bphi$ or $\bpsi$ exceed a given threshold $\tau$, these two vectors will be normalized with the excessive part `absorbed' in $\boldsymbol{\alpha}$ and $\boldsymbol{\beta}$: 
    \begin{equation*}
         \boldsymbol{\alpha} \gets \boldsymbol{\alpha} + \ln(\bphi^{(\ell)}),\quad
         \boldsymbol{\beta}  \gets \boldsymbol{\beta}  + \ln(\bpsi^{(\ell)}), \quad
         \bphi^{(\ell)} \gets \mathbf{1}_N, \quad \bpsi^{(\ell)} \gets \mathbf{1}_N.
     \end{equation*}
     Correspondingly, the matrix $K$ needs to be rescaled as $K \gets diag(e^{\boldsymbol{\alpha}/\veps})\times K\times diag(e^{\boldsymbol{\beta}/\veps}).$
\end{remark}
	
\section{The FS-1 Algorithm}
\label{sec:3}
    The key to the Sinkhorn algorithm is to iteratively update $\bphi$ and $\bpsi$ through equation \eqref{update_sinkhorn}. By introducing the notation $\lambda = e^{-h/\veps}$, the matrix multiplication vector operation is written as
	\begin{equation}\label{eqn:NBS_summation}
		K \bpsi^{(\ell)} =
		\begin{pmatrix}
			\psi^{(\ell)}_1 &+& \lambda\psi^{(\ell)}_2 &+& \lambda^2\psi^{(\ell)}_3 & \cdots &+& \lambda^{N-1}\psi^{(\ell)}_N\\
			\lambda\psi^{(\ell)}_1 &+& \psi^{(\ell)}_2 &+& \lambda\psi^{(\ell)}_3 &\cdots  &+& \lambda^{N-2}\psi^{(\ell)}_N\\
			\lambda^2\psi^{(\ell)}_1&+& \lambda\psi^{(\ell)}_2&+&  \psi^{(\ell)}_3  &\cdots &+&\lambda^{N-3}\psi^{(\ell)}_N \\
			\vdots &\vdots & \vdots &\vdots & \vdots & \ddots &\vdots & \vdots\\
			\lambda^{N-1}\psi^{(\ell)}_1&+& \lambda^{N-2}\psi^{(\ell)}_2&+& \lambda^{N-3}\psi^{(\ell)}_3 & \cdots& +&\psi^{(\ell)}_N
		\end{pmatrix}.
	\end{equation}
	We separate the summation of row $k$ to the lower triangular part $p_k$ and the strictly upper triangular part $q_k$. Then updating $\bphi^{(\ell+1)}$ is formulated as
	\begin{equation*}\label{meshpoly}
		\phi^{(\ell+1)}_k = u_k/(p_k+q_k), \quad
		p_k = \sum_{i=1}^{k} \psi^{(\ell)}_i \lambda^{k-i}, \quad
		q_k = \sum_{i=k+1}^{N} \psi^{(\ell)}_i \lambda^{i-k}, \quad k=1,2,\cdots,N.
	\end{equation*}
	Instead of directly calculating $p_k$ and $q_k$, we use the recursive computation  given by
	\begin{equation}\label{recursion}
		\begin{aligned}
			&p_1 = \psi^{(\ell)}_1, \quad 
			p_{k+1} = \lambda p_k +\psi^{(\ell)}_{k+1}, \quad k=1,2,\cdots,N-1, \\
			&q_N = 0, \quad
			q_{k} = \lambda(q_{k+1}+\psi^{(\ell)}_{k+1}), \quad k=N-1, N-2, \cdots,1.
		\end{aligned}
	\end{equation}
	This is the Qin Jiushao or Horner method for efficient polynomial evaluation \cite{needham1974, horner1819}.
	Thus we develop the FS-1 algorithm with linear computational complexity, which only takes $2(N-1)$ times additions and multiplications for the matrix multiplication operation. The pseudo-code is presented in Algorithm \ref{FS-1}.
	
	\begin{algorithm}[h]
		\caption{FS-1 Algorithm}
		\label{FS-1}
		\hspace*{0.02in} {\bf Input:} $\bu$ and $\bv$ of size $(N,1)$,\ $h$,\ $\veps$,\ $\mathrm{tol}$,\ $\mathrm{itr\_max}$ \\
		\hspace*{0.02in} {\bf Output:} $W_\veps(\bu,\bv)$
		\begin{algorithmic}[1]
			\State $\lambda=e^{-h/\veps}$; $\bphi^{(0)},\bpsi^{(0)} = \frac{1}{N}\mathbf{1}_N$; $\bp^{(0)},\br^{(0)},\bq^{(0)},\bs^{(0)}=\mathbf{0}_N$; $\ \ell=0$
			\While{$\brac{\ell\;<\;\mathrm{itr\_max}}$ and $\brac{\sum_{i=1}^N\abs{\psi_i^{(\ell)}(\sum_{j=1}^N\phi_j^{(\ell)}\lambda^{\abs{i-j}})-v_i}\;>\;\mathrm{tol}}$}
			\State $r^{(\ell)}_{1} =\phi^{(\ell)}_1$, $s^{(\ell)}_{N} = 0$			
			\For{$i=1\;:\;N-1$}
			\State $r^{(\ell)}_{i+1} = \lambda r^{(\ell)}_{i}+\phi^{(\ell)}_{i+1}$
			\State $s^{(\ell)}_{N-i} = \lambda(s^{(\ell)}_{N-i+1}+\phi^{(\ell)}_{N-i+1})$
			\EndFor
			\State $\bpsi^{(\ell+1)}=\bv\oslash\brac{\br^{(\ell)}+\bs^{(\ell)}}$
			\State $p^{(\ell)}_{1} =\psi^{(\ell+1)}_1$, $q^{(\ell)}_{N} = 0$
			\For{$i=1\;:\;N-1$}
			\State $p^{(\ell)}_{i+1} = \lambda p^{(\ell)}_{i}+\psi^{(\ell+1)}_{i+1}$
			\State $q^{(\ell)}_{N-i} = \lambda(q^{(\ell)}_{N-i+1}+\psi^{(\ell+1)}_{N-i+1})$
			\EndFor
			\State $\bphi^{(\ell+1)}=\bu\oslash\brac{\bp^{(\ell)}+\bq^{(\ell)}}$
			\State $\ell=\ell+1$
			\EndWhile{}			
			\noindent\Return $\sum_{i,j=1}^N \phi^{(\ell)}_i \psi^{(\ell)}_j \lambda^{\abs{i-j}}\abs{i-j}h$
		\end{algorithmic}
	\end{algorithm}	

	It is well-known that the Sinkhorn algorithm has stability issues~\cite{chizat2018scaling} due to the division and the multiplication of the small parameter $\lambda^k$. Next, we need to discuss the stability of the FS-1 algorithm to ensure that its stability is not worse than the Sinkhorn algorithm.
	
	Consider the matrix multiplication in \eqref{eqn:NBS_summation}, 
	\begin{equation*}
	    S^{\alpha}=K\bpsi^{\alpha}, \quad
	    S^{\alpha}=(S^{\alpha}_1,\;S^{\alpha}_2, \cdots, S^{\alpha}_N)^T, \quad
	    \bpsi^{\alpha}=(\psi^{\alpha}_1,\;\psi^{\alpha}_2, \cdots, \psi^{\alpha}_N)^T, \quad \alpha=1,2.
	\end{equation*}
	assume that
	\begin{equation*}
	    \abs{\psi^1_k-\psi^2_k}\le \delta, \quad k=1,2,\cdots,N,
	\end{equation*}
	thus 
    \begin{multline*}
        \abs{S_k^1-S_k^2}=\abs{\lambda^{k-1}\brac{\psi_1^1-\psi_1^2}\!+\!\lambda^{k-2}\brac{\psi_2^1-\psi_2^2}\!+\!\cdots\!+\!
        \brac{\psi_k^1-\psi_k^2}\!+\!\cdots\!+\!\lambda^{N-k}\brac{\psi_N^1-\psi_N^2}}, \\
        \le \abs{\lambda^{k-1}+\lambda^{k-2}+\cdots+\lambda+1+\lambda+\cdots+\lambda^{N-k}}\delta
        =\brac{\frac{2-\brac{\lambda^k+\lambda^{N-k+1}}}{1-\lambda}-1}\delta.
    \end{multline*}
    On the other hand, consider the successive computation in  \eqref{recursion}, an easy induction gives
    \begin{equation*}
        \abs{p_k^1-p_k^2}\le \lambda\abs{p_{k-1}^1-p_{k-1}^2}+\abs{\psi_k^1-\psi_k^2}
        \le \cdots \le \brac{1+\lambda+\cdots\lambda^{k-1}}\delta
        =\frac{1-\lambda^k}{1-\lambda}\delta,
    \end{equation*}
    and
    \begin{equation*}
        \abs{q_k^1-q_k^2}\le \lambda\abs{q_{k\!+\!1}^1-q_{k\!+\!1}^2}+\abs{\psi_{k\!+\!1}^1-\psi_{k\!+\!1}^2} \le \cdots \le \brac{\lambda+\cdots+\lambda^{N\!-\!k\!-\!1}}\delta
        \le \brac{\frac{1-\lambda^{N\!-\!k}}{1-\lambda}-1}\delta.
    \end{equation*}
    Thus, we have
    \begin{equation*}
        \abs{(p_k^1+q_k^1)-(p_k^2+q_k^2)}
        \le \brac{\frac{2-\brac{\lambda^k+\lambda^{N-k+1}}}{1-\lambda}-1}\delta,
    \end{equation*}
    that is, the FS-1 algorithm and the Sinkhorn algorithm has the same stability.

\medskip

\begin{remark}	
    Similarly, the log-domain stabilization \cite{chizat2018scaling} technique can also be aggregated into the FS-1 algorithm. First, we need to `absorb' the excessive part of vectors $\bphi$ and $\bpsi$ into $\boldsymbol{\alpha}$ and $\boldsymbol{\beta}$:
    \begin{equation*}
         \boldsymbol{\alpha} \gets \boldsymbol{\alpha} + \ln(\bphi^{(l)}),\quad
         \boldsymbol{\beta}  \gets \boldsymbol{\beta}  + \ln(\bpsi^{(l)}), \quad
         \bphi^{(l)} \gets \mathbf{1}_N, \quad \bpsi^{(l)} \gets \mathbf{1}_N.
     \end{equation*}
     Next, we have to replace lines 5-6 and 10-11 in the Algorithm \ref{FS-1} as 
    \begin{align*}\label{differ_alg}
	   5: \quad r^{(\ell)}_{i+1} &= \lambda e^{(\beta_{i+1}-\beta_{i})/\veps}r^{(\ell)}_{i}+e^{(\alpha_{i+1}+\beta_{i+1})/\veps}\phi^{(\ell)}_{i+1}, \\
	   6: \quad s^{(\ell)}_{N-i} &= \lambda e^{(\beta_{N-i}-\beta_{N-i+1})/\veps}\brac{s^{(\ell)}_{N-i+1}+e^{(\alpha_{N-i+1}+\beta_{N-i+1})/\veps}\phi^{(\ell)}_{N-i+1}}, \\
	   10: \quad p^{(\ell)}_{i+1} &= \lambda e^{(\alpha_{i+1}-\alpha_{i})/\veps}p^{(\ell)}_{i}+e^{(\alpha_{i+1}+\beta_{i+1})/\veps}\psi^{(\ell+1)}_{i+1}, \\
	   11: \quad q^{(\ell)}_{N-i} &= \lambda e^{(\alpha_{N-i}-\alpha_{N-i+1})/\veps}\brac{q^{(\ell)}_{N-i+1}+e^{(\alpha_{N-i+1}+\beta_{N-i+1})/\veps}\psi^{(\ell+1)}_{N-i+1}}.
	\end{align*}
	This leads to the stabilized FS-1 algorithm with $O(N)$ complexity.
\end{remark}

\medskip

 \begin{remark} \label{rmk:stab_ana}
     In practice, the direct calculation of $\lambda^m=e^{-mh/\veps}$ and the multiplication of $\lambda^m$ may be troublesome since $\lambda^m$ could be very small. However, the FS-1 algorithm gets around this problem by stepwise multiplication of $\lambda$.
 \end{remark}

\section{Extension to high dimension}
\label{sec:4}
In this section, we illustrate how the FS-1 algorithm generalizes to higher dimensions using the two-dimensional case as an example. First, consider two discretized probabilistic distributions 
\begin{align*}
    \bu &= \brac{u_{11},u_{21},\cdots,u_{N1},u_{12},\cdots,u_{i_1j_1}, \cdots, u_{NM}},\\
    \bv &= \brac{v_{11},v_{21},\cdots,v_{N1},v_{12},\cdots, v_{i_2j_2}, \cdots, v_{NM}},
\end{align*}
on a uniform 2D mesh of size $N\times M$ with a vertical spacing of $h_1$ and a horizontal spacing of $h_2$, the entropy regularized 2D Wasserstein-1 metric can be discretized as the optimal value of the following minimizing problem:
\begin{equation}\label{eqn:2d disc entropy problem}
	W_\veps(\bu,\bv) = \!\inf_{\Gamma_{\!i_1j_1i_2j_2\!}} \sum_{\!i_1,i_2=1\!}^N\sum_{\!j_1,j_2=1\!}^M \brac{\Gamma_{\!i_1j_1i_2j_2\!}\brac{\!\abs{i_1-i_2}h_1\!+\!\abs{j_1-j_2}h_2\!}+\veps\Gamma_{\!i_1j_1i_2j_2\!}\ln\big( \Gamma_{\!i_1j_1i_2j_2\!} \big)}\!,
 \end{equation} 
where $\Gamma_{i_1j_1i_2j_2}$ satisfies
\begin{equation*}
    \sum_{i_2=1}^N\sum_{j_2=1}^M\Gamma_{i_1j_1i_2j_2} = u_{i_1j_1},\quad\sum_{i_1=1}^N\sum_{j_1=1}^M\Gamma_{i_1j_1i_2j_2} = v_{i_2j_2},\quad\Gamma_{i_1j_1i_2j_2}\geq0.
\end{equation*}
Same as in the 1D case, the problem \eqref{eqn:2d disc entropy problem} can be solved by the Sinkhorn iteration. If $\bu$ and $\bv$ are flattened into 1D vectors in column-major order, the corresponding kernel matrix is written as
\begin{equation*}
\renewcommand\arraystretch{2}
    	K = \left(\begin{array}{c|c|c|c|c}
		K_0 & \lambda_2 K_0 & \lambda_2^2K_0 & \cdots & \lambda_2^{M-1}K_0 \\
		\hline
		\lambda_2 K_0& K_0 & \lambda_2 K_0 &\cdots & \lambda_2^{M-2}K_0 \\
		\hline
		\vdots &\vdots &\vdots  &\ddots & \vdots \\
		\hline
		\lambda_2^{M-1}K_0\; & \lambda_2^{M-2}K_0\; & \lambda_2^{M-3}K_0\; & \cdots & K_0
	\end{array}\right),
\end{equation*}
where the sub-matrix
\begin{equation*}
    K_0=\begin{pmatrix}
	1 & \lambda_1 &\cdots &\lambda_1^{N-1} \\
	\lambda_1 & 1 &\cdots &\lambda_1^{N-2} \\
	\vdots	& \vdots &\ddots &\vdots	\\
	\lambda_1^{N-1} & \lambda_1^{N-2} & \cdots & 1
\end{pmatrix},
\end{equation*}
and
\begin{equation*}
	\lambda_1=e^{-h_1/\veps}, \quad \lambda_2 = e^{-h_2/\veps}.
\end{equation*}

Obviously, the cost of direct matrix-vector multiplication in Sinkhorn is $O(N^2M^2)$. By using the FS-1's trick twice, we expect the computational cost can be significantly reduced, the key idea is as follows. Let
\begin{equation*}
	\bpsi^{(\ell)}_i=\brac{\psi_{1i}^{(\ell)},\psi_{2i}^{(\ell)},\cdots,\psi_{Ni}^{(\ell)}},\quad i=1,2,\cdots,M,
\end{equation*}
 the matrix-vector multiplication $K\bpsi^{(\ell)}$ is written as
\begin{equation}\label{eqn:NBS_summation 2d}
	K \bpsi^{(\ell)} =
	\begin{pmatrix}
		K_0\bpsi^{(\ell)}_1 &+& \lambda_2 K_0\bpsi^{(\ell)}_2 &+& \lambda_2^2K_0\bpsi^{(\ell)}_3 & \cdots &+& \lambda_2^{M-1}K_0\bpsi^{(\ell)}_M\\
		\lambda_2 K_0\bpsi^{(\ell)}_1 &+& K_0\bpsi^{(\ell)}_2 &+& \lambda_2 K_0\bpsi^{(\ell)}_3 &\cdots  &+& \lambda_2^{M-2}K_0\bpsi^{(\ell)}_M\\
		\lambda_2^2K_0\bpsi^{(\ell)}_1&+& \lambda_2 K_0\bpsi^{(\ell)}_2&+&  K_0\bpsi^{(\ell)}_3  &\cdots &+&\lambda_2^{M-3}K_0\bpsi^{(\ell)}_M \\
		\vdots &\vdots & \vdots &\vdots & \vdots & \ddots &\vdots & \vdots\\
		\lambda_2^{M-1}K_0\bpsi^{(\ell)}_1&+& \lambda_2^{M-2}K_0\bpsi^{(\ell)}_2&+& \lambda_2^{M-3}K_0\bpsi^{(\ell)}_3 & \cdots& +&K_0\bpsi^{(\ell)}_M
	\end{pmatrix}.
\end{equation}
We separate the summation of row $k$ to the lower triangular part $\bp_k$ and the strictly upper triangular part $\bq_k$. Then updating $\bphi^{(\ell+1)}$ is formulated as
\begin{equation*}\label{meshpoly 2d}
	\bphi^{(\ell+1)}_k = \bu_k\oslash(\bp_k+\bq_k), \quad
	\bp_k = \sum_{i=1}^{k} K_0\bpsi^{(\ell)}_i \lambda_2^{k-i}, \quad
	\bq_k = \sum_{i=k+1}^{M} K_0\bpsi^{(\ell)}_i \lambda_2^{i-k}, \quad k=1,\cdots,M.
\end{equation*}
Instead of directly calculating $\bp_k$ and $\bq_k$, a successive computation is used
\begin{equation}\label{2d recursion}
	\begin{aligned}
		&\bp_1 = K_0\bpsi^{(\ell)}_1, \quad 
		\bp_{k+1} = \lambda_2 \bp_k + K_0\bpsi^{(\ell)}_{k+1}, \quad k=1,2,\cdots,M-1, \\
		&\bq_M = \mathbf{0}_N, \quad
		\bq_{k} = \lambda_2(\bq_{k+1}+K_0\bpsi^{(\ell)}_{k+1}), \quad k=M-1,M-2,\cdots,1.
	\end{aligned}
\end{equation}
Thus, the total cost of matrix-vector multiplication of our FS-1 algorithm for 2D Wasserstein-1 metric is reduced to $O(NM)$. The pseudo-code is presented in Algorithm \ref{FS-12d}. 

This idea can be easily extended to high-dimensional cases, and we will not repeat it here. An argument similar to the one used in Section 3 shows that the FS-1 algorithm and the Sinkhorn algorithm in the 2D case still has the same stability. Thus, we shall omit the discussions.
	
	\begin{algorithm}[h]
		\caption{2D FS-1 Algorithm}
		\label{FS-12d}
		\hspace*{0.02in} {\bf Input:} $\bu$ and $\bv$ of size $(N\times M,1)$,\ $h_1$,\ $h_2$,\ $\veps$,\ $\mathrm{tol}$,\ $\mathrm{itr\_max}$ \\
		\hspace*{0.02in} {\bf Output:} $W_\veps(\bu,\bv)$
		\begin{algorithmic}[1]
			\State $\lambda_1=e^{-h_1/\veps};\lambda_2=e^{-h_2/\veps}$; $\bphi^{(0)},\bpsi^{(0)} = \frac{1}{NM}\mathbf{1}_{NM}$; $\bp^{(0)},\br^{(0)},\bq^{(0)},\bs^{(0)}=\mathbf{0}_{NM}$; $\ell=0$
			\While{$\brac{\ell\;<\;\mathrm{itr\_max}}$ and $\brac{\sum_{i=1}^M\Vert\bpsi_i\odot(\sum_{j=1}^M \lambda_2^{\abs{i-j}}K_0\bphi_j^{(\ell)})-\bv_i\Vert_1\;>\;\mathrm{tol}}$}
			\State $\br^{(\ell)}_{1} =K_0\bphi^{(\ell)}_1$ (using Algorithm \ref{FS-1} with $h=h_1$,\;$\veps=\veps$), $\bs^{(\ell)}_{M} = \mathbf{0}_N$			
			\For{$i=1\;:\;M-1$}
			\State $\br^{(\ell)}_{i+1} = \lambda_2 \br^{(\ell)}_{i}+K_0\bphi^{(\ell)}_{i+1}$
			\State $\bs^{(\ell)}_{M-i} = \lambda_2(\bs^{(\ell)}_{M-i+1}+K_0\bphi^{(\ell)}_{M-i+1})$
			\EndFor
			\State $\bpsi^{(\ell+1)}=\bv\oslash\brac{\br^{(\ell)}+\bs^{(\ell)}}$
			
			\State $\bp^{(\ell)}_{1} =K_0\bpsi^{(\ell+1)}_1$ (using Algorithm \ref{FS-1} with $h=h_1$,\;$\veps=\veps$), $\bq^{(\ell)}_{M} = \mathbf{0}_N$
			\For{$i=1\;:\;M-1$}
			\State $\bp^{(\ell)}_{i+1} = \lambda_2 \bp^{(\ell)}_{i}+K_0\bpsi^{(\ell+1)}_{i+1}$
			\State $\bq^{(\ell)}_{M-i} = \lambda_2(\bq^{(\ell)}_{M-i+1}+K_0\bpsi^{(\ell+1)}_{M-i+1})$
			\EndFor
			\State $\bphi^{(\ell+1)}=\bu\oslash\brac{\bp^{(\ell)}+\bq^{(\ell)}}$
			\State $\ell=\ell+1$
			\EndWhile{}			
			\noindent\Return $\sum\limits_{i_1,i_2=1}^N \sum\limits_{j_1,j_2=1}^M \phi^{(\ell)}_{i_1,j_1} \psi^{(\ell)}_{i_2,j_2} \lambda_1^{\abs{i_1-i_2}}\lambda_2^{\abs{j_1-j_2}}(\abs{i_1-i_2}h_1+\abs{j_1-j_2}h_2)$
		\end{algorithmic}
	\end{algorithm}	

\section{Numerical Experiments}
\label{sec:5}
In this section, we carry out four numerical experiments to evaluate the FS-1 algorithm. The first two examples show the performance of the FS-1 algorithm in 1D cases. Specifically, we consider the comparison of two 1D random distributions and the comparison of two Ricker wavelets arising from seismology \cite{chen2018quadratic}. The last two examples show the performance of the FS-1 algorithm in 2D cases. Specifically, we consider the comparison of two 2D random distributions and the comparison of two images arising from the image matching problem \cite{ILSVRC15}. The marginal error $\Vert \text{diag}(\bpsi)\times K^T\times \bphi-\bv\Vert_1$ is chosen as the termination condition \cite{altschuler2017near,scetbon2021low}. All the experiments are conducted on a platform with 128G RAM, and one Intel(R) Xeon(R) Gold 5117 CPU @2.00GHz with 14 cores. 

\subsection{1D random distributions}
\label{subsec:1drandom}
We consider the Wasserstein-1 metric between two 1D random distributions on the interval $[-3,3]$. There are uniform grid points
\begin{equation*}
	x_i=(i-1)\Delta x-3, \quad \Delta x=\frac{6}{N-1},\quad i=1,\;2,\cdots,N.
\end{equation*}
Correspondingly, we consider the two random vectors on the grid points
\begin{equation*}
	\bu=(u_1,\;u_2,\cdots,u_N), \quad \bv=(v_1,\;v_2,\cdots,v_N),
\end{equation*}
where $u_i$ and $v_j$ are both uniformly distributed on $[0,1]$. We would like to compare the performance and computational cost on computing the Wasserstein-1 metric $W_{\veps}\brac{\frac{\bu}{\norm{\bu}},\frac{\bv}{\norm{\bv}}}$ using the Sinkhorn algorithm and the FS-1 algorithm. We tested 100 random experiments, and each experiment was performed for 1000 iterations.

\begin{figure}
    \centering
    \includegraphics[width=\linewidth]{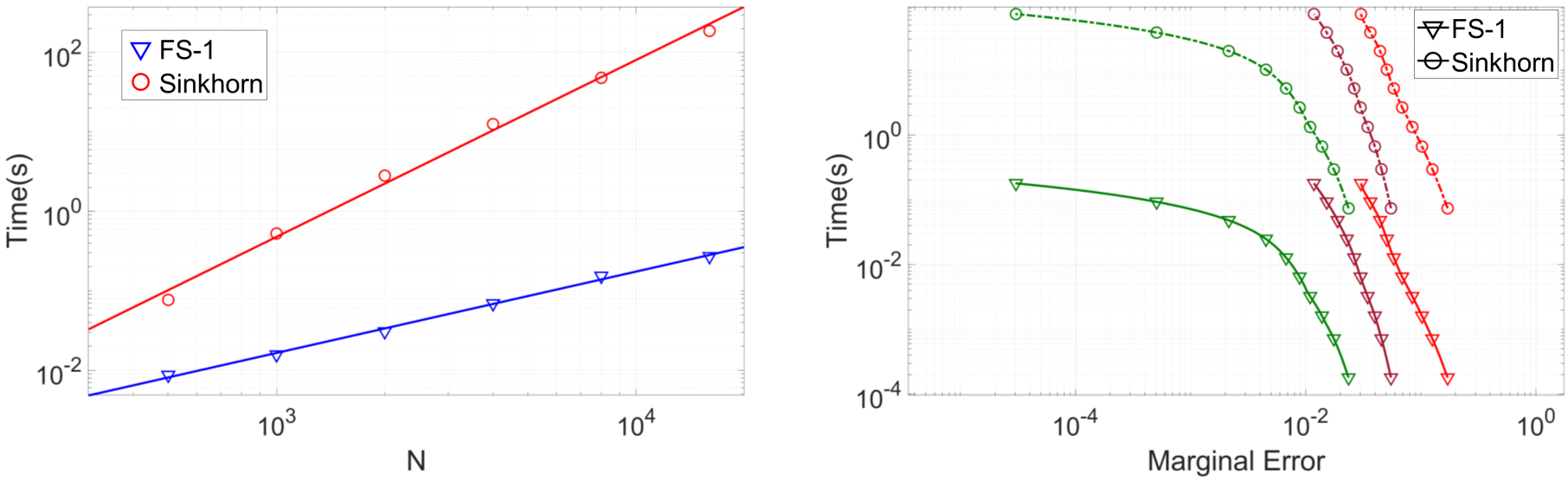}
    \caption{The 1D random distribution problem. Left: The comparison of computational time between the FS-1 algorithm and the Sinkhorn algorithm with different numbers of grid pointss $N$. Right: The computational time required to reach the corresponding marginal error under different regularization parameter $\veps=0.1$(Green), $\veps=0.01$(Purple), and $\veps=0.001$(Red). }
    \label{fig:1drandom}
\end{figure}

In Table \ref{randomOTexp}, we output the averaged computational time of two different algorithms. We can see that the FS-1 algorithm has an overwhelming advantage in computational speed. Moreover, the transport plans obtained by the two algorithms are almost identical. To further study the efficiency advantage of our FS-1 algorithm, we present the computational time of the two algorithms in different cases. By data fitting, we can see the empirical complexity of the FS-1 algorithm is $O(N^{1.02})$, while that of the Sinkhorn algorithm is $O(N^{2.22})$, see Figure \ref{fig:1drandom} (Left) for illustration. In Figure \ref{fig:1drandom} (Right), we discuss the computational time required to reach the corresponding marginal error under different regularization parameters $\veps$ for the random distribution with dimension $N=10000$. Obviously, the FS-1 algorithm is about two orders of magnitude faster than the Sinkhorn algorithm.

\begin{table}[h]
\centering
	\begin{tabular}{cclcc}
		\toprule
		\multirow{2}{*} N &
		\multicolumn{2}{c}{Computational time (s)} & \multirow{2}{*}{Speed-up ratio} & \multirow{2}{*}{$\norm{P_{NB}-P}_F$} \\
		
		\cline{2-3} & FS-1 & Sinkhorn \\
		
		 \midrule
		500   & $8.70\times10^{-3}$ & $7.68\times10^{-2}$ & $8.83\times10^{0}$  & $6.54\times10^{-15}$   \\ 
		2000  & $4.25\times10^{-2}$ & $2.81\times10^{0}$  & $6.61\times10^{1}$ &$4.98\times10^{-18}$   \\ 
		8000  & $1.53\times10^{-1}$ & $4.80\times10^{1}$ & $3.14\times10^{2}$ &$3.92\times10^{-18}$   \\ 
		\bottomrule
	\end{tabular}
	\caption{The 1D random distribution problem. The comparison between the Sinkhorn algorithm and the FS-1 algorithm with the different number of grid points $N$. The  regularization  parameter $\veps=0.001$. Columns 2-4 are the averaged computational time of the two algorithms and the speed-up ratio of the FS-1 algorithm. Column 5 is the Frobenius norm of the difference between the transport plan computed by the two algorithms.}
	\label{randomOTexp}
\end{table}

\subsection{Ricker wavelet}

Next, we consider the computation of the Wasserstein-1 metric between the Ricker wavelet
\begin{equation*}
	R(t)=A(1-2\pi^2f_0^2t^2)e^{-\pi^2f_0^2t^2},
\end{equation*}
and its translation $R(t-s)$. The Ricker wavelet is commonly used to model source time function in seismology \cite{chen2018quadratic}. Here $f_0$ is the dominant frequency, and $A$ denotes the wave amplitude. For simplicity, we set
\begin{equation*}
	f_0=1, \quad A=1.
\end{equation*}
Since the Ricker wavelet is not always positive over the entire time duration, we will  square and normalize it for the comparisons of the Wasserstein-1. In \cite{lzy2021}, a new normalization method with better convexity is given as follows:
\begin{equation} \label{ricker}
	D(f,g)=W_{1,\veps}\brac{\frac{\frac{f^2}{\norm{f^2}}+\delta}{1+L\delta},
		\frac{\frac{g^2}{\norm{g^2}}+\delta}{1+L\delta}},
\end{equation}
where $\delta$ is a sufficiently small parameter to improve numerical stability, and $L$ is a given parameter to guarantee that the two functions being compared are normalized.

\begin{figure}
    \centering
    \includegraphics[width=\linewidth]{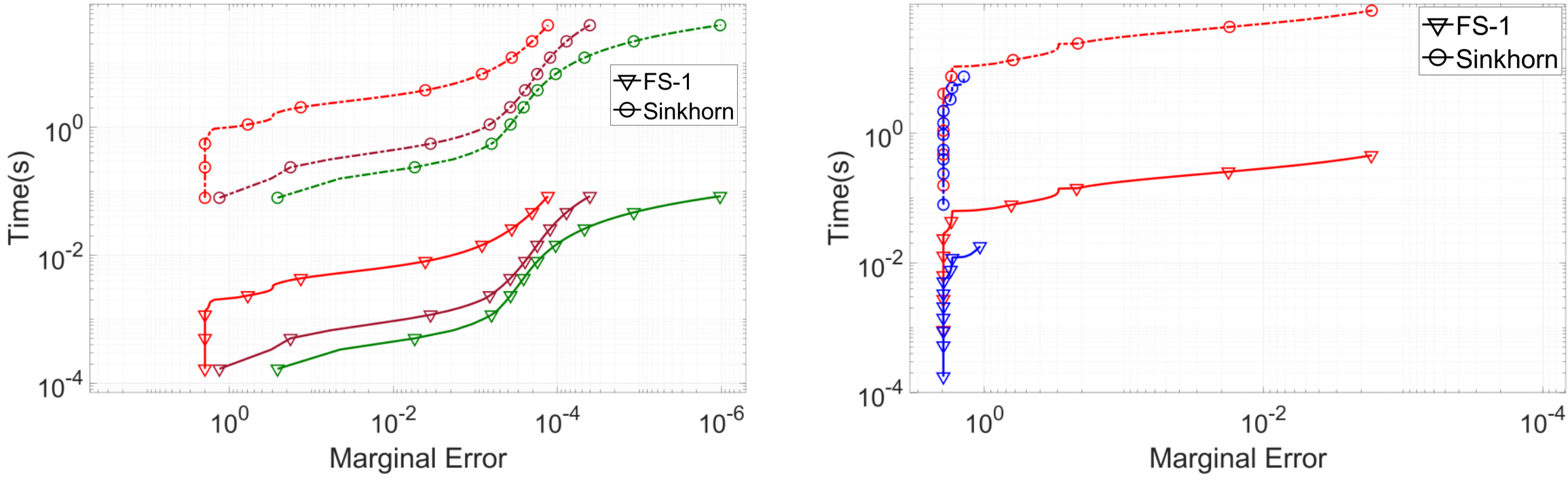}
    \caption{The Ricker wavelet problem. Left: The computational time required to reach the corresponding marginal error under different regularization parameter $\veps=0.1$(Green), $\veps=0.05$(Purple), and $\veps=0.01$(Red). Right: The comparison between the Sinkhorn-type algorithms with (Red) and without (Blue) the log-domain stabilization for $\veps=0.001$.}
    \label{fig:1dricker}
\end{figure}

Below, we randomly select the translation parameter $s=-1.2032$. This parameter can ensure that the two Ricker wavelets $R(t)$ and $R(t-s)$ are sufficiently far apart. And the small parameter $\delta=10^{-3}$. We repeated the experiment 100 times, and each experiment was performed for 500 iterations. In Table \ref{rickerOTexp}, we output the averaged computational time of two different algorithms. We also present the computational time required to reach the corresponding marginal error under different regularization parameters $\veps$ for the distribution with dimension $N=10000$, see Figure \ref{fig:1dricker} (Left) for illustration, from which, we can draw the same conclusions as those in Subsection \ref{subsec:1drandom}.

In Figure \ref{fig:1dricker} (Right), we also discuss the impact of the log-domain stabilization technique for $\veps=0.001$. Without the technique, the Sinkhorn algorithm terminates abnormally at the 97th iteration and the FS-1 algorithm terminates abnormally at the 104th iteration. By introducing the log-domain stabilization technique, neither algorithm is terminated abnormally. Moreover, the stabilized FS-1 algorithm still maintains a significant efficiency advantage over the stabilized Sinkhorn algorithm.

\begin{table}[h]
\centering
	\begin{tabular}{cclcc}
		\toprule
		\multirow{2}{*}N &
		\multicolumn{2}{c}{Computational time (s)} & \multirow{2}{*}{Speed-up ratio} & \multirow{2}{*}{$||P_{NB}-P||_F$} \\
		
		\cline{2-3} & FS-1 &	Sinkhorn \\
		
		 \midrule
		500   & $4.90\times10^{-3}$ & $2.30\times10^{-2}$ & $4.69\times10^{0}$  & $5.67\times10^{-16}$   \\ 
		2000  & $1.52\times10^{-2}$ & $1.36\times10^{0}$  & $8.93\times10^{1}$ &$1.81\times10^{-17}$   \\ 
		8000  & $5.96\times10^{-2}$ & $2.81\times10^{1}$ & $4.72\times10^{2}$ &$1.22\times10^{-16}$   \\ 
		\bottomrule
	\end{tabular}
	\caption{The Ricker wavelet problem. The comparison between the Sinkhorn algorithm and the FS-1 algorithm with the different number of grid points $N$. The regularization  parameter $\veps=0.01$. Columns 2-4 are the averaged computational time of the two algorithms and the speed-up ratio of the FS-1 algorithm. Column 5 is the Frobenius norm of the difference between the transport plan computed by the two algorithms.}
	\label{rickerOTexp}
\end{table}

\subsection{2D Random distributions}
\label{subsec:2drandom}
\begin{figure}
    \centering
    \includegraphics[width=\linewidth]{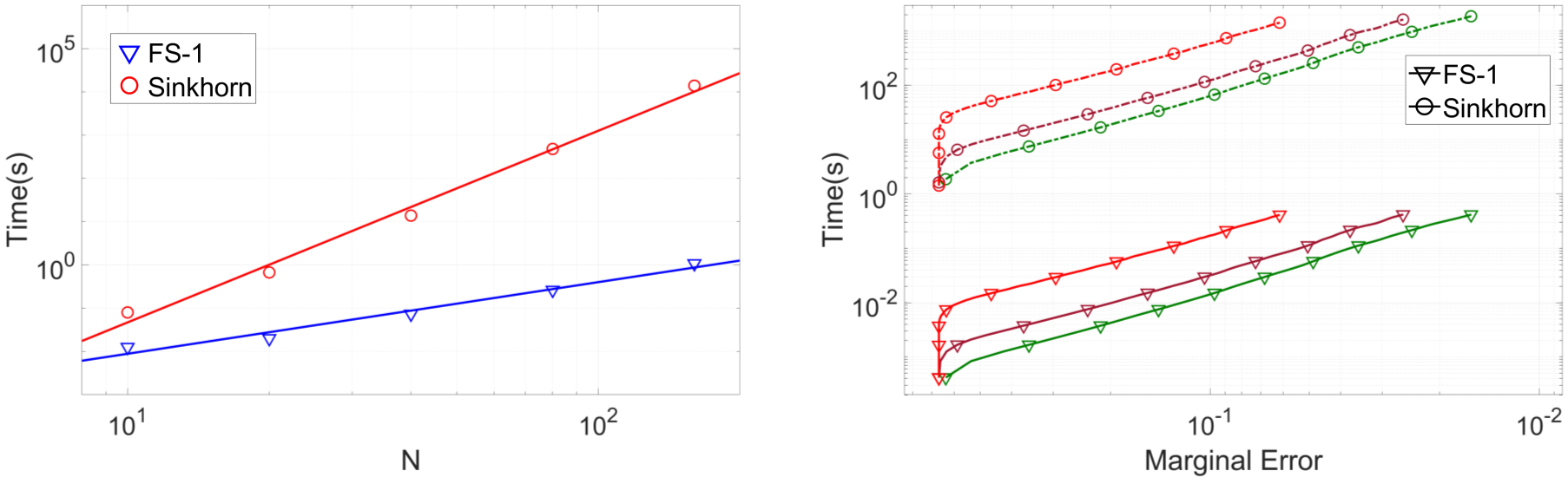}
    \caption{The 2D random distribution problem. Left: The comparison of computational time between the FS-1 algorithm and the Sinkhorn algorithm with different numbers of grid nodes $N$. Right: The computational time required to reach the corresponding marginal error under different regularization parameter $\veps=0.1$(Green), $\veps=0.05$(Purple), and $\veps=0.01$(Red).}
    \label{fig:2drandom}
\end{figure}

Next, we discuss the performance of the FS-1 algorithm in two dimension. This Subsection generalizes the 1D random distributions in Subsection \ref{subsec:1drandom} to 2D random distributions. The basic settings are almost the same as in Subsection \ref{subsec:1drandom}.  And we need to compute the Wasserstein-1 metric between two $N\times N$ dimensional random vectors, where $N$ is the number of grid points at each dimension. Without loss of generality, we set $h_1=h_2=1$. We also tested 100 random experiments, and each experiment was performed for 1000 iterations. The averaged computational time for different numbers of nodes $N\times N$ of the two algorithms are output in Table \ref{2drandomOTexp} and Figure \ref{fig:2drandom} (Left). By data fitting, we can see the empirical complexity of the FS-1 algorithm is $O(N^{1.66})$ , while that of the Sinkhorn algorithm is $O(N^{4.44})$. These results are even better than  the expected $O(N^2)$ complexity for the FS-1 algorithm. This again shows the big efficiency advantage of the FS-1 algorithm compared to the Sinkhorn algorithm. In Figure \ref{fig:2drandom} (Right), we also present the computational time required to reach the corresponding marginal error under different regularization parameters $\veps$ for the random distribution with dimension $100\times100$. Obviously, the FS-1 algorithm still maintains the efficiency  advantage for more than two orders of magnitude.

\begin{table}[h]
\centering
	\begin{tabular}{ccccc}
		\toprule
		\multirow{2}{*}{$N\times N$} &
		\multicolumn{2}{c}{Computational time (s)} & \multirow{2}{*}{Speed-up ratio} & \multirow{2}{*}{$||P_{NB}-P||_F$} \\
		
		\cline{2-3} & FS-1 &	Sinkhorn \\
		
		 \midrule
		10$\times$10   & $1.27\times10^{-2}$ & $8.02\times10^{-2}$ & $6.34\times10^{0}$  & $1.20\times10^{-17}$   \\ 
		20$\times$20  & $2.01\times10^{-2}$ & $6.73\times10^{-1}$  & $3.34\times10^{1}$ &$5.96\times10^{-18}$   \\ 
		40$\times$40  & $7.38\times10^{-2}$ & $1.38\times10^{1}$ & $1.87\times10^{2}$ &$3.00\times10^{-18}$   \\ 
		80$\times$80  & $2.64\times10^{-1}$ & $4.78\times10^{2}$ & $1.81\times10^{3}$ &$1.55\times10^{-18}$   \\
		160$\times$160  & $1.08\times10^{0}$ & $1.38\times10^{4}$ & $1.28\times10^{4}$ &$7.68\times10^{-19}$   \\
		\bottomrule
	\end{tabular}
	\caption{The 2D random distribution problem. The comparison between the Sinkhorn algorithm and the FS-1 algorithm with the different number of grid points $N\times N$. The regularization  parameter $\veps=0.01$. Columns 2-4 are the averaged computational time of the two algorithms and the speed-up ratio of the FS-1 algorithm. Column 5 is the Frobenius norm of the difference between the transport plan computed by the two algorithms.}
	\label{2drandomOTexp}
\end{table}

\subsection{Image matching problem}

An important application of the optimal transport in 2D is to match images. OT plays a fundamental role of related tasks including density regularization~\cite{burger2012regularized}, image registration~\cite{haker2004optimal}, and optical flow~\cite{clarysse2010optimal}. The complexity advantage of the FS-1 algorithm compared to the Sinkhorn algorithm can further enable practical applications of optimal transport in high-resolution images. 

\begin{figure}[h!]
\centering
	\includegraphics[width=0.3\linewidth,height=0.25\linewidth]{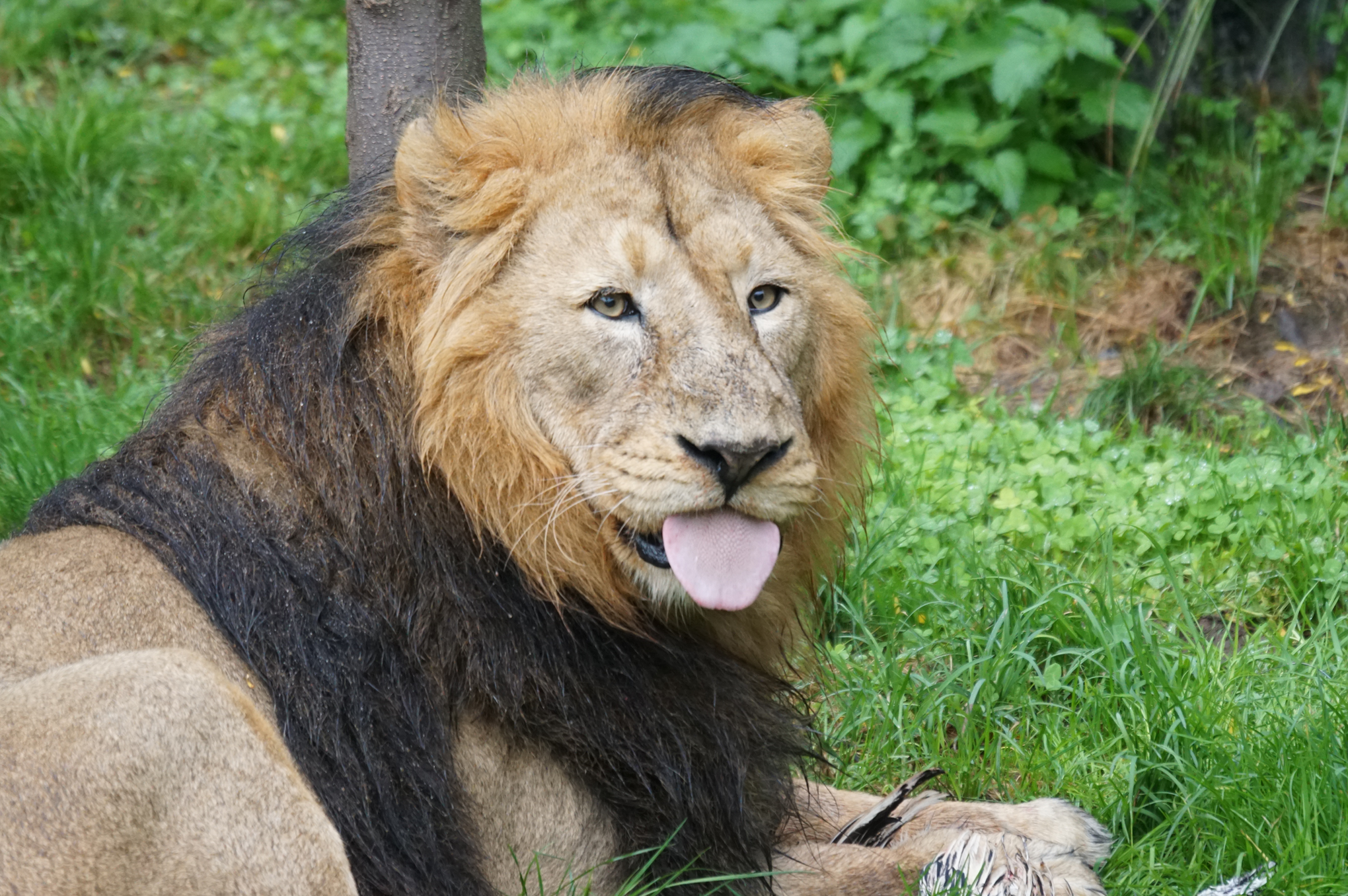}
	\includegraphics[width=0.3\linewidth,height=0.25\linewidth]{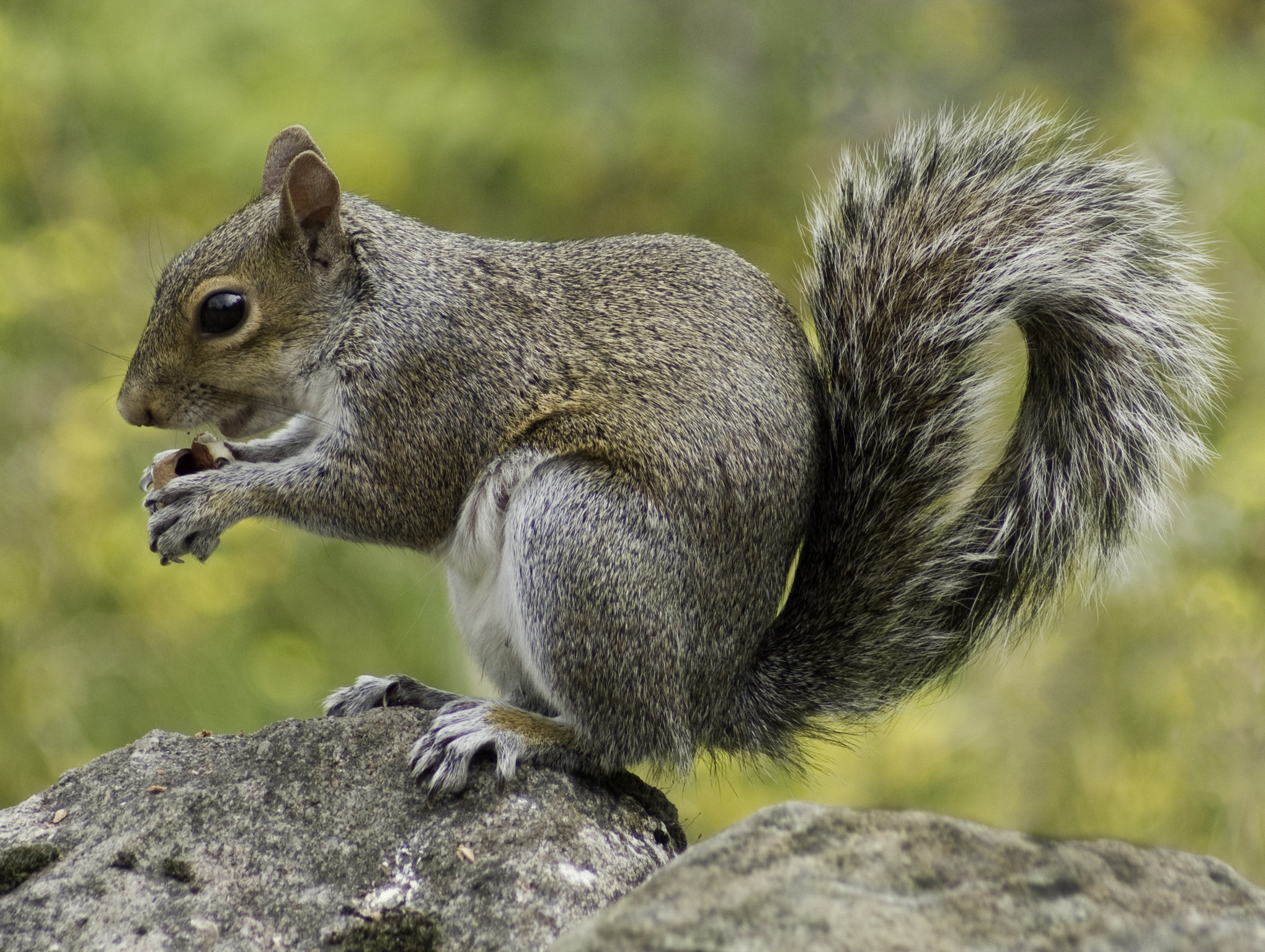}
\caption{The image matching problem. Illustration of images.}
\label{ilsvrc}
\end{figure}

Here we consider the image matching experiment. We randomly select two images, see Figure \ref{ilsvrc} for illustration, from the DIV2K dataset\cite{Agustsson_2017_CVPR_Workshops}, where the images have 2K pixels for at least one of the axes (vertical or horizontal). Considering that the resolutions of the images are different, we first sample them to the same scale $N\times N$. Without loss of generality, we set $h_1=h_2=1$. In addition, to facilitate the optimal transport comparison, we convert them to grayscale images and normalize them using  formula \eqref{ricker} with $\delta=10^{-7}$. Again, we repeated the experiment $100$ times and each experiment was performed for $1000$ iterations.

\begin{figure}
    \centering
    \includegraphics[width=1\linewidth]{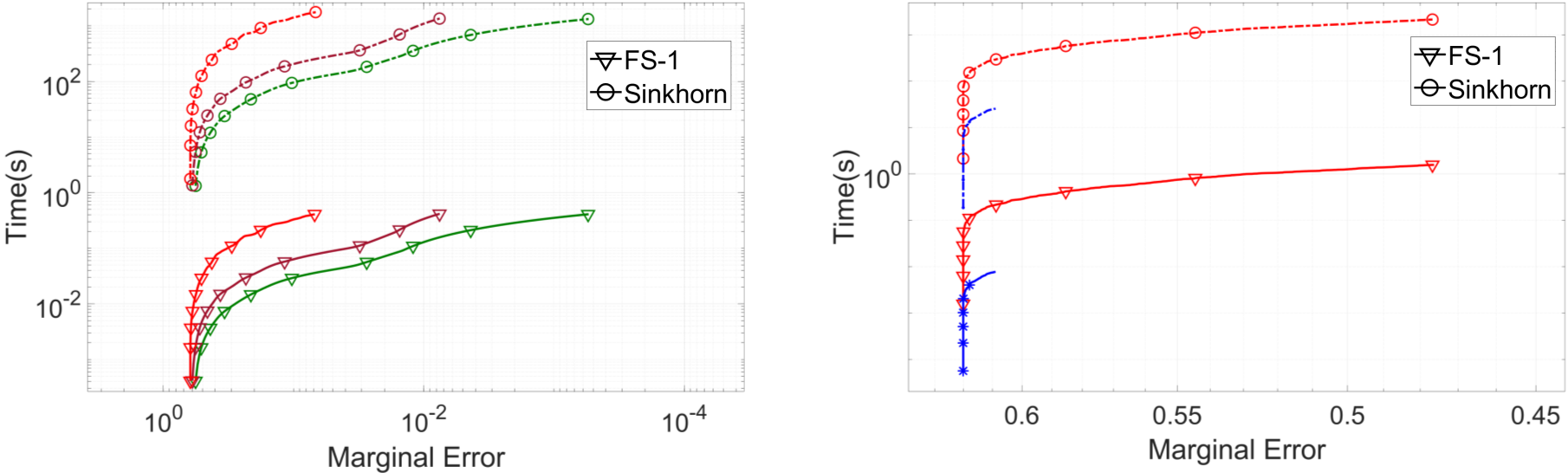}
    \caption{The image matching problem. Left: The computational time required to reach the corresponding marginal error under different regularization parameter $\veps=1$(Green), $\veps=0.5$(Purple), and $\veps=0.1$(Red). Right: The comparison between the Sinkhorn-type algorithms with (Red) and without (Blue) the log-domain stabilization for $\veps=0.01$.}
    \label{fig:2dfig}
\end{figure}

In Table \ref{2dfigtable}, we output the averaged computational time of two different algorithms. It should be emphasized that the computational time of the Sinkhorn algorithm is too long for large-scale images, so there is no result. However, the FS-1 algorithm can still  handle this simply. We also present the computational time required to reach the corresponding marginal error under different regularization parameters $\veps$. Here, the image size is $100\times100$, see Figure \ref{fig:2dfig} (Left) for illustration,  from which  we can draw the same conclusions as those in Subsection \ref{subsec:2drandom}.

In Figure \ref{fig:2dfig} (Right), we also discuss the impact of the log-domain stabilization technique for $\veps=0.01$. Without the technique, the Sinkhorn algorithm and the FS-1 algorithm both terminate abnormally at the 138th iteration. By introducing the log-domain stabilization technique, neither algorithm is terminated abnormally. Moreover, the stabilized FS-1 algorithm still maintains a significant efficiency advantage over the stabilized Sinkhorn algorithm.



\begin{table}[h]
\centering
	\begin{tabular}{ccccc}
		\toprule
		\multirow{2}{*}{$N\times N$} &
		\multicolumn{2}{c}{Computational time (s)} & \multirow{2}{*}{Speed-up ratio} & \multirow{2}{*}{$\norm{P_{NB}-P}_F$} \\
		
		\cline{2-3} & FS-1 &	Sinkhorn \\
		
		 \midrule
		100$\times$100  & $4.10\times10^{-1}$ & $1.32\times10^{3}$ & $3.23\times10^{3}$  & $2.28\times10^{-17}$   \\ 
		200$\times$200  & $1.72\times10^{0}$ & $3.08\times10^{4}$  & $1.79\times10^{5}$ &$9.52\times10^{-18}$   \\ 
		400$\times$400  & $7.23\times10^{0}$ & $-$ & $-$ &$-$   \\ 
		800$\times$800  & $3.20\times10^{1}$ & $-$ & $-$ &$-$   \\
		\bottomrule
	\end{tabular}
	\caption{The  image  matching  problem. The  comparison  between  the  Sinkhorn  algorithm  and the FS-1 algorithm with the different total number of grid nodes $N\times N$.  The regularization parameter $\veps=1$. Columns 2-4 are the averaged computational time of the two algorithms and the speed-up ratio of the FS-1 algorithm.  Column 5 is the Frobenius norm of the difference between the transport plan computed by the two algorithms.}
	\label{2dfigtable}
\end{table}

\section{Conclusion}
\label{sec:6}
In this paper, we propose an efficient (abbreviated as FS-1) algorithm to compute the Wasserstein-1 metric with linear computational cost per iteration. This method is developed by discovering the natural structure of Sinkhorn's kernel, which allows  a matrix-vector multiplication to be carried out exactly with $O(n)$ cost for each iteration by using Qin Jiushao's or Horner's method for efficient polynomial evaluation.  Moreover, the FS-1 algorithm can also be  adapted to the widely used  log-domain stabilization technique. As shown by numerous experiments, the FS-1 algorithm achieves a huge speed advantage without losing accuracy. 

Finally, this paper mainly considers the acceleration of the Sinkhorn algorithm in the matrix-vector multiplication. It is well known that the number of iterations of the Sinkhorn algorithm will significantly increase with the increase of numerical accuracy, which leads to slow convergence. In~\cite{xie2020fast}, the Inexact Proximal point method for the Optimal Transport problem (IPOT) was proposed for this problem. We believe that FS-1 and IPOT can be effectively combined to develop a new algorithm for solving the Optimal Transport problem with fast convergence and low complexity. We are currently investigating this important extension and hope to report the progress in a future paper.



\section*{Acknowledgements}

This work was supported by National Natural Science Foundation of China Grant Nos.11871297 and 12031013, Tsinghua University Initiative Scientific Research Program, and Shanghai Municipal Science and Technology Major Project 2021SHZDZX0102.

\bibliographystyle{siam}
\bibliography{ref.bib}

\end{sloppypar}
\end{document}